\documentclass{amsart}
\usepackage{amssymb,amsmath}

\newtheorem{theorem}{Theorem}
\newcommand{\bt}{\begin{theorem}}
\newcommand{\et}{\end{theorem}}

\newtheorem{lemma}{Lemma}
\newcommand{\bl}{\begin{lemma}}
\newcommand{\el}{\end{lemma}}

\newtheorem{corollary}{Corollary}
\newcommand{\bc}{\begin{corollary}}
\newcommand{\ec}{\end{corollary}}

\newtheorem*{conjecture}{Conjecture}
\newcommand{\bconj}{\begin{conjecture}}
\newcommand{\econj}{\end{conjecture}}

\newtheorem{example}{Example}
\newcommand{\bex}{\begin{example}}
\newcommand{\eex}{\end{example}}

\newtheorem{problem}{Problem} 
\newcommand{\bp}{\begin{problem}}
\newcommand{\ep}{\end{problem}}

\newcommand{\beq}{\begin{equation}}
\newcommand{\eeq}{\end{equation}}
\newcommand{\benum}{\begin{enumerate}}
\newcommand{\eenum}{\end{enumerate}}
\newcommand{\ba}{\begin{array}}
\newcommand{\ea}{\end{array}}
\newcommand{\Z}{\ensuremath{\mathbf Z}}

\newcommand{\PP}{\ensuremath{\mathbf P}}

\newcommand{\F}{\ensuremath{\mathbf F }}

\DeclareMathOperator{\card}{card}

\begin{document}
\date{\today}
\title[Heights in finite projective space]{Heights in finite projective space, and a problem on directed graphs} 
\author{Melvyn B. Nathanson}
\address{Department of Mathematics\\
Lehman College (CUNY)\\
Bronx, New York 10468}
\email{melvyn.nathanson@lehman.cuny.edu}
\thanks{The work of M.B.N. was supported in part by grants from the NSA Mathematical Sciences Program and the PSC-CUNY Research Award Program.}
\author{Blair D. Sullivan}
\address{Department of Mathematics\\
Princeton University\\
Princeton, New Jersey 08544}
\email{bdowling@princeton.edu}
\thanks{The work of B.D. S. was supported in part by a Department of Homeland Security Dissertation Grant.}

\begin{abstract}
Let $\F_p = \Z/p\Z$.  The \emph{height} of a point $\mathbf{a}=(a_1,\ldots, a_d) \in \F_p^d$ is $h_p(\mathbf{a}) = \min \left\{\sum_{i=1}^d (ka_i \mod p) : k=1,\ldots,p-1\right\}.$ Explicit formulas and estimates are obtained for the values of the height function in the case $d=2,$ and these results are applied to the problem of determining the minimum number of edges the must be deleted from a finite directed graph so that the resulting subgraph is acyclic.

\end{abstract}

\date{\today}

\maketitle

\section{Heights in finite projective space}
Let $F$ be a field and let $F^{\ast} = F\setminus\{0\}$.  For 
$d \geq 2$, we define an equivalence relation on the set of nonzero 
$d$-tuples 
$F^{d}\setminus\{(0,\ldots,0)\}$ 
as follows: 
$(a_1,\ldots,a_{d}) \sim (b_1,\ldots,b_d)$ 
if there exists $k \in F^{\ast}$ such that 
$(b_1,\ldots,b_{d})=(ka_1,\ldots,ka_{d})$.
 We denote the equivalence class of $(a_1,\ldots,a_{d})$
by $\langle a_1,\ldots,a_{d}\rangle$.
The set of equivalence classes is called the $(d-1)$-dimensional projective space over the field $F$, and denoted $\PP^{d-1}(F)$. 

We consider projective space over the finite field $\F_p = \Z/p\Z.$  For every $x \in \F_p$, we denote by  $x \mod p$  the least nonnegative  integer in the congruence class $x$.  We define the \emph{height} of the point  $\mathbf{a} = \langle a_1,\ldots,a_{d}\rangle \in \PP^{d-1}(\F_p)$ by
\[
h_p(\mathbf{a}) = \min \left\{\sum_{i=1}^d (ka_i \mod p) : k=1,\ldots,p-1\right\}.
\]
For every nonempty set $\mathcal{A} \subseteq \PP^{d-1}(\F_p)$, we define
\[
H_p(\mathcal{A}) = \{ h_p(\mathbf{a}) : \mathbf{a} \in \mathcal{A} \}.
\]
Then $H_p(\mathcal{A})$ is a set of positive integers.

For $\mathbf{a} = \langle a_1,\ldots,a_{d}\rangle \in \PP^{d-1}(\F_p)$, let $d^{\ast}(\mathbf{a})$ denote the number of nonzero components of $\mathbf{a},$ that is, the number of $a_i \neq 0.$   The function $d^{\ast}(\mathbf{a})$ is well-defined, that is, independent of the representative of the equivalence class of $\mathbf{a}.$  For $\mathcal{A} \subseteq \PP^{d-1}(\F_p)$, we define
\[
d^{\ast}(\mathcal{A}) = \max\{ d^{\ast}(\mathbf{a}) : \mathbf{a} \in \mathcal{A} \}.
\]
Then $h_p(\mathbf{a}) \leq d^{\ast}(\mathbf{a})(p-1)$ for all $\mathbf{a} \in \PP^{d-1}(\F_p)$.  We can reduce this upper bound by a simple averaging argument.

For every real number $t$, let $[t]$ denote the greatest integer not exceeding $t$.

\bl  \label{CSS:lemma:upperbound}
For every point $\mathbf{a} \in \PP^{d-1}(\F_p)$,
\[
h_p(\mathbf{a}) \leq \left[ \frac{d^{\ast}(\mathbf{a})p}{2} \right].
\]
\el

\begin{proof}
If $a \in \F_p^{\ast}$, then $\{ka \mod p : k=1,\ldots,p-1\} = \{1,\ldots,p-1\}$ and so 
\[
\sum_{k=1}^{p-1} (ka \mod p) = \sum_{k=1}^{p-1} k = \frac{p(p-1)}{2}.
\]
It follows that for every $\mathbf{a} = \langle a_1,\ldots,a_{d}\rangle \in \PP^{d-1}(\F_p)$, we have
\[
\sum_{k=1}^{p-1} \sum_{i=1}^d (ka_i \mod p)
 = \sum_{i=1}^d \sum_{k=1}^{p-1} (ka_i \mod p)  = \frac{d^{\ast}(\mathbf{a}) p(p-1)}{2}.
\]
Since the minimum of a set of numbers does not exceed the average of the set, we have
\[
h_p(\mathbf{a}) \leq \frac{1}{p-1}\sum_{k=1}^{p-1} \sum_{i=1}^d (ka_i \mod p)= \frac{d^{\ast}(\mathbf{a})p}{2}.
\]
The Lemma follows from the fact that the heights are positive integers.
\end{proof}

\bl  \label{CSS:lemma:HA}
For every odd prime $p$ and  $d \geq 2,$
\begin{align*}
\max\left(  H_p( \PP^{d-1}(\F_p) ) \right)  = \frac{dp}{2} & \qquad\text{if $d$ is even} \\
\frac{(d-1)p}{2}+1 \leq \max\left( H_p( \PP^{d-1}(\F_p) ) \right) \leq \frac{dp-1}{2}
& \qquad\text{if $d$ is odd.}
\end{align*}
\el

\begin{proof}
If $2r \leq d$ and $a_1,\ldots,a_r,a_{2r+1},\ldots,a_d$ are nonzero elements of the field $\F_p$, then the point 
\[
\mathbf{a} = \langle a_1,a_2,\ldots,a_{r},-a_1,-a_2\ldots,-a_{r},a_{2r+1},\ldots,a_d\rangle,
\]
satisfies $d^{\ast}(\mathbf{a}) = d$ and
\[
\begin{split}
\sum_{i=1}^d (ka_i  &\mod p) \\
 &= \sum_{i=1}^{r} \left((ka_i \mod p) + (-ka_i \mod p) \right) + \sum_{i=2r+1}^{d} (ka_i \mod p) \\
& \geq rp + d-2r
\end{split}
\]
for all $k=1,\ldots,p-1$.  If $d-2r \leq p-1,$ then we can choose distinct elements $a_{2r+1},\ldots,a_d$ and 
\[
\sum_{i=1}^d (ka_i  \mod p)  \geq rp+ \frac{(d-2r)(d-2r+1)}{2}.
\]
Applying Lemma~\ref{CSS:lemma:upperbound} and the inequality with $r=[d/2]$, we obtain
$h_p(\mathbf{a}) = dp/2$
if $d$ is even and 
\[
\frac{(d-1)p}{2}+1 \leq h_p(\mathbf{a}) \leq \frac{dp-1}{2}
\]
if $d$ is odd.  This completes the proof.
\end{proof}

\section{Heights on the finite projective line}
The projective line $\PP^1(\F_p)$ consists of all equivalence classes of pairs $(a_1,a_2)$, where $a_1,a_2 \in \F_p$ and $a_1$ and $a_2$ are not both 0.
If $a_1=0$, then $\langle 0,a_2\rangle  = \langle 0,1\rangle $ and $h_p(\langle 0,1\rangle )=1.$
If $a_2=0$, then $\langle a_1,0\rangle  = \langle 1,0\rangle $ and $h_p(\langle 1,0\rangle )=1.$
If $a_1\neq 0$ and $a_2 \neq 0$, then $\langle a_1,a_2\rangle  = \langle 1,a_1^{-1}a_2\rangle $.  
Thus, for all $\mathbf{a} \in \PP^1(\F_p)$, if $\mathbf{a} \neq \langle 1,0\rangle $ and $\mathbf{a} \neq \langle 0,1\rangle $, then $\mathbf{a} = \langle 1,a\rangle $ for some $a \in \F_p^{\ast},$ and $h_p(\langle 1,a\rangle ) \geq 2.$

\bl     \label{CSS:lemma:1a}
Let $p$ be an odd prime and $a \in \F_p^{\ast}.$  Then
\benum
\item[(i)]
$h_p(\langle 1,a\rangle )  \leq 1+ (a \mod p)$  for all $a$,
\item[(ii)]
$h_p(\langle 1,a\rangle ) = 1+ (a \mod p)$     if  $a\mod p < \sqrt{p}$,
\item[(iii)]
$h_p(\langle 1,a\rangle ) = 2$ if and only if $a=1+p\Z$,
\item[(iv)]
$h_p(\langle 1,a\rangle )  = 3$ if and only if $a = 2+p\Z$ or $a = (p+1)/2+ p\Z$,
\item[(v)]
$h_p(\langle 1,a\rangle ) = p$    if and only if $a = p-1 + p\Z$,
\item[(vi)]
Let $a = p-b + p\Z$ for $1 \leq b \leq p-1.$  Then
$h_p( \langle 1,a\rangle )  \leq \left(p+(b-1)^2\right)/ b.$
\eenum
\el 

\begin{proof}
For all $a \in \F_p^{\ast}$ and $k \in \{1,\ldots,p-1\}$ we have $ka \mod p \in \{1,\ldots,p-1\}$, and so
\[
h_p(\langle 1,a\rangle )  = \min\{k + (ka\mod p): k=1,\ldots,p-1\}  \leq 1 +( a \mod p) .
\]
Note that $ka\mod p \leq k(a\mod p)$ for all $k \geq 1.$  If $k \geq a \mod p$, then $k + (ka\mod p) \geq  (a\mod p) +1.$  
If $1 \leq k \leq (a \mod p)-1$ and $(a\mod p) < \sqrt{p}$, then 
\[
ka \mod p  \leq k(a \mod p) \leq ((a \mod p)-1)(a \mod p) \leq (a \mod p)^2 < p
\]
It follows that $ka \mod p = k(a \mod p)$ and
\[
k+(ka\mod p)  = k + k(a \mod p) \geq 1 + (a\mod p).
\]
and so $h_p(\langle 1,a\rangle ) = 1+(a\mod p).$  This proves~(i) and~(ii).

We have $k + (ka\mod p) = 2$ if and only if $k=1$ and $ka\mod p = a\mod p = 1$, that is, $a= 1+p\Z.$  Similarly, $k + (ka\mod p) = 3$ if and only if either $k=1$ and $ka\mod p = a\mod p = 2$, or $k=2$ and $ka\mod p = 2a\mod p = 1$.  In the first case, $a= 2+p\Z$ and, in the second case, $a= (p+1)/2 + p\Z.$    This proves~(iii) and~(iv).

If $a=-1+p\Z$, then $k + (ka \mod p) = k + (p-k) = p$ for all $k=1,\ldots,p-1$ and so $h_p(1,a) = p.$  Conversely, if $h_p(1,a) = p,$ then $k + (ka\mod p)= p$ for some $k$, and so $ka \mod p = -k \mod p$ and $a=-1+p\Z.$    This proves~(v). 

Finally, to prove~(vi), we let $p=qb+r,$ where $q= [p/b]$ and $1 \leq r \leq p-1.$  Then
\[
qa  = \left[ \frac{p}{b} \right] (p-b) + p\Z  = p- \left[ \frac{p}{b} \right]b + p\Z = r+p\Z
\]
and so $qa \mod p = r.$
Therefore,
\[
h_p(\langle 1,a)\rangle ) \leq q + r = \frac{p+r(b-1)}{b}  \leq \frac{p + (b-1)^2}{b}.
\]
This completes the proof.
\end{proof}

\bt    \label{CSS:theorem:(p-1)/2}
Let $p$ be an odd prime and $a \in \F_p.$  Then
$h_p(\langle 1,a\rangle )  = (p+1)/2$ if and only if $a = (p-1)/2+p\Z$ or $a = p-2 +p\Z$.
If 
\[
a \notin \{ (p-1)/2+ p\Z,p-2+p\Z,p-1+p\Z \},
\]
then
\[
h_p(\langle 1,a\rangle )  \leq \frac{p-1}{2} .
\]
\et

\begin{proof}
The Theorem is true for $p = 3,5,$ and 7, so we can assume that $p \geq 11.$

Let $a = p-2+p\Z.$  If $1 \leq k \leq (p-1)/2,$ then 
\[
k + (ka \mod p) = k + (p-2k) = p-k \geq \frac{p+1}{2}
\]
and $k + (ka \mod p) =  (p+1)/2$ when $k = (p-1)/2.$
If $k \geq (p+1)/2,$ then $k + (ka \mod p) \geq (p+3)/2.$  Therefore, $h_p(\langle 1,a)\rangle ) = (p+1)/2.$

Let $a = (p-1)/2 + p\Z$.  If $j=1,\ldots, (p-1)/2$ and $k=2j,$ then
\[
k + (ka\mod p) = 2j + (j(p-1)\mod p)  = 2j + (p-j) = p+j \geq p+1.
\]
If $k=2j-1,$ then
\begin{align*}
k + (ka\mod p) & = (2j-1) + \left( \frac{  (2j-1)(p-1)  }{2} \mod p \right)\\
& = (2j-1) + \left( \frac{p+1}{2} - j \right) \\
& = \frac{p+2j-1}{2} \geq  \frac{p+1}{2}.
\end{align*}
Since $1 + (a\mod p) = (p+1)/2,$ it follows that $h_p(\langle 1,a\rangle ) =  (p+1)/2.$

If $a \in \F_p^{\ast}$ and $a \mod p \in \{0,1,2,\ldots, (p-3)/2\},$ then 
\[
h_p(\langle 1,a\rangle )  \leq 1 + (a  \mod p)\leq \frac{p-1}{2}
\]
by Lemma~\ref{CSS:lemma:1a} (i).
If $a \in \F_p^{\ast}$ and $a \mod p = (p+1)/2,$ then 
$h_p(\langle 1,a\rangle ) = 3< (p+1)/2$ by Lemma~\ref{CSS:lemma:1a} (iv).

Let $a \in \F_p^{\ast}$ and $(p+3)/2 \leq a \mod p \leq p-3.$   
There is an integer $b$ such that
\[
3 \leq b \leq \frac{p-3}{2} \qquad \text{and} \qquad a = p-b + p\Z.
\]
By Lemma~\ref{CSS:lemma:1a} (vi) we have 
$h_p( \langle 1,a\rangle )  \leq \left(p+({b}-1)^2\right)/ {b},$ and so 
$h_p( \langle 1,a\rangle )  \leq (p-1)/2$ if 
\[
2b+1 + \frac{4}{b-2} \leq p.
\]
If $4 \leq b \leq (p-3)/2,$ then 
\[
2b+1 + \frac{4}{b-2} \leq 2b+3 \leq p.
\]
If $b = 3,$ then $h_p(\langle 1,a\rangle ) = h_p(\langle 1,p-3\rangle ) \leq (p-1)/2$ since 
\[
2b+1 + \frac{4}{b-2} = 11 \leq p.
\]
This completes the proof.
\end{proof}

\section*{Table of heights for primes $11 \leq p \leq 29$}

\begin{center} 
\begin{tabular}{|c|c|c|c|c|c|}
\hline
prime $p$ & $ a \mod p $ & $h_p(\langle 1,a\rangle )$ & 
prime $p$ & $ a \mod p $ & $h_p(\langle 1,a\rangle )$ \\ \hline\hline
11 & 2 & 3 & 23&  2&  3\\
     & 3 & 4 & & 3 & 4 \\
     & 4 & 4 & & 4 &  5\\
     & 5 &  6 & & 5 &  6\\
     & 6 & 3 & &  6&  5\\
     & 7 & 5 & &  7&  8\\
     & 8 & 5 & &  8&  4\\
     & 9 & 6 & &  9&  7\\  \cline{1-3}
13 & 2 & 3 & &  10&  8\\
     & 3 & 4 & &  11&  12\\
     & 4 & 5 & &  12&  3\\
     & 5 &  5 & &  13& 5 \\ 
     & 6 &  7& &  14&  6\\ 
     & 7 &  3& &  15&  9\\
     & 8 &  5 & &  16&  5\\
     & 9 &  4 & &  17&  8\\
     & 10 & 5 & &  18&  7\\
     & 11 & 7 & &  19&  8\\  \cline{1-3}
17 & 2 & 3 & &  20&  9\\
     & 3 & 4 & &  21&  12\\  \cline{4-6}
     & 4 &  5 & 29&  2& 3 \\
     & 5 &   6 & & 3 & 4 \\ 
     & 6 &  4& &  4&  5\\ 
     & 7 &  6 & & 5 &  6\\ 
     & 8 &   9 & & 6 &  6\\ 
     & 9 &  3 & &  7&  8\\
     & 10 &  5 & & 8 &  7\\
     & 11 &  7 & & 9 &  10\\
     & 12 &  5 & & 10 &  4\\
     & 13 & 5 & &  11&  7\\
     & 14 & 7& &  12 &  7\\
     & 15 & 9 & &  13 &  10\\  \cline{1-3}
19 & 2 & 3 & &  14 &  15\\
     & 3 & 4 & &  15&  3\\
     & 4 &  5& &  16&  5\\
     & 5 &   5 & &  17&  7\\ 
     & 6 & 7 & &  18 &  8\\ 
     & 7 &  5 & &  19 & 11 \\ 
     & 8 &   7 & &  20&  5\\ 
     & 9 &   10 & &  21 & 8 \\ 
     & 10 &   3& &  22 & 5 \\
     & 11 &  5 & &  23& 9 \\
     & 12 &   7& &  24 & 9 \\
     & 13 &  4& &  25& 8 \\
     & 14 & 7& &  26& 11 \\
     & 15 &  7& &  27& 15 \\
     & 16 & 7&  &  &  \\
     & 17 & 10& &  &  \\
     \hline
\end{tabular}
\end{center}

\section{Problems on heights}
\bp
Let $d \geq 2$ and $\mathbf{a} = \langle a_1,\ldots,a_{d}\rangle \in \PP^{d-1}(\F_p)$.  Is there a simple formula  to compute $h_p(\mathbf{a})$?  
Is there a simple formula  to estimate $h_p(\mathbf{a})$?  This is not known even for the projective line $d=2$.  
\ep

\bp
By Theorem~\ref{CSS:theorem:(p-1)/2} and Lemma~\ref{CSS:lemma:1a}, we have 
\[
 H_p( \PP^{1}(\F_p) ) \bigcap \left( \frac{p+1}{2} , p\right) = \emptyset.
 \]
For which positive integers $r$ does there exist a number $c_r$ such that 
\[
 H_p( \PP^{1}(\F_p) ) \bigcap \left( \frac{p}{r+1}+c_r , \frac{p}{r}-c_r\right) = \emptyset
 \]
for all sufficiently large $p$?
\ep

\bp
Is there an upper bound for the heights of points in the projective plane $\PP^2(\F_p)$ analogous to the upper bound in Theorem~\ref{CSS:theorem:(p-1)/2} for the projective line?
\ep

\bp  \label{HeightsLine:problem}
The following problem arises in graph theory.  Let  $k \geq 2$ and let $\mathcal{A} \subseteq \PP^{d-1}(\F_p)$ be a nonempty subset of projective space such that 
\benum
\item
If  $\mathbf{a} = \langle a_1,\ldots,a_{d}\rangle \in \PP^{d-1}(\F_p)$, then the coordinates $a_i$ are pairwise distinct.
\item
For $\ell = 1,\ldots, k,$ none of the equations
\[
x_1+x_2  + \cdots + x_{\ell}  = 0 
\]
has a solution with $x_1,\ldots, x_k\in \{a_1,a_2,\ldots,a_d\}.$  (These conditions are homogeneous and  independent of the representative of the equivalence class of $\mathbf{a}$.)  
\eenum
Find an upper bound for $H_p(\mathcal{A}).$
\ep

\bp
Find a good definition of the height of a point in the projective space $\PP^{d-1}(\F_q)$ over any finite field $\F_q.$
\ep

\section{Cayley graphs with vertex set $\F_p$} 
Let $G = (V,E)$ be a directed graph with vertex set $V$ and edge set $E \subseteq V\times V.$  
A \emph{directed path} of length $n$ in $G$ is a sequence of vertices $v_{i_0},v_{i_1},v_{i_2},\ldots, v_{i_n}$ such that $(v_{i_j},v_{i_{j+1}})$ is an edge for $j=0,1,\ldots, n-1$.  
A \emph{directed cycle} of length $n$ in $G$ is a directed path  $v_{i_0},v_{i_1},v_{i_2},\ldots, v_{i_n}$ such that $v_{i_n}=v_{i_0}$.  A loop is a cycle of length 1, a digon is a cycle of length 2, and a triangle is a cycle of length 3.  A 3-free or triangle-free graph is a graph with no loop, digon, or triangle.  The graph  $G = (V,E)$ is called  \emph{directed acyclic} if it has no directed cycle.

The outdegree of the vertex $v$ is the number of edges of the form $(v,v')$ for some vertex $v '.$  The pigeonhole principle implies that in a finite directed graph, if the outdegree of every vertex is at least 1, then the graph contains a cycle.    Thus, every finite directed acyclic graph contains at least one vertex with outdegree 0.

\bt  \label{CSS:theorem:permutation}
Let $\{k_0,k_1,\ldots, k_{m-1}\}$ be a set of $m$ distinct integers, and let $G$ be a finite directed graph with vertex set $V = \{v_{k_0}, v_{k_1},\ldots, v_{k_{m-1}}\}$.  The graph $G$ is directed acyclic if and only if there is a one-to-one map $\sigma: \{0,1,\ldots, m-1\} \rightarrow \{k_0,k_1,\ldots, k_{m-1}\}$ such that, if $(v_{\sigma(i)},v_{\sigma(j)})$ is an edge of the graph, then $i < j.$  If 
$\{k_0,k_1,\ldots, k_{m-1}\} = \{0,1,\ldots, m-1\},$  then $G$ is directed acyclic if and only if there is a permutation $\sigma$ of $\{0,1,\ldots, m-1\}$ such that $r < s$ for every edge $(v_{\sigma(r)},v_{\sigma(s)})$ of the graph.
\et

\begin{proof}
Let $\sigma: \{0,1,\ldots, m-1\} \rightarrow \{k_0,k_1,\ldots, k_{m-1}\}$ be a one-to-one map such that, if $(v_{\sigma(i)},v_{\sigma(j)})$ is an edge of the graph, then $i<j.$  If 
$v_{\sigma(i_0)},v_{\sigma(i_1)},v_{\sigma(i_2)},\ldots, v_{\sigma(i_n)}$ is a path in $G$, then $i_0 < i_1 < i_2 < \cdots < i_n$ and so $i_n \neq i_0,$ that is, $v_\sigma(i_n) \neq v_\sigma(i_0)$, and so no path in $G$ is a cyclic.

To prove the converse, we use induction on $m$.  The Lemma holds for $m=1$ and $m=2.$.   Assume that $m\geq 2$ and that the Lemma is true for every finite acyclic graph with  $m$ vertices.  If $G$ is an acyclic directed graph with $m+1$ vertices $\{v_{k_0}, v_{k_1},\ldots, v_{k_{m}}\}$, then there exists a vertex $v_{k_r}$ with outdegree 0.   Consider the induced subgraph $G'$ of $G$ on the vertex set $\{ v_{k_0},v_{k_1},\ldots, v_{k_{r-1}}, v_{k_{r+1}},\ldots, v_{k_m} \}$.  By the induction hypothesis, there is a one-to-one map $\sigma'$ from $\{0,1,\ldots, m-1\}$ into $\{k_0,k_1,\ldots,k_{r-1},k_{r+1},\ldots, k_m\}$ such that if $(v_{\sigma'(i)},v_{\sigma'(j)})$ is an edge of the graph $G'$, then $\sigma'(i) < \sigma'(j).$   Extend this map to a function $\sigma$ of $\{0,1,\ldots, m\}$ by defining $\sigma(i) = \sigma'(i)$ for $i=0,1,\ldots,m-1$ and $\sigma(m) = k_r.$  Since $v_{k_r} = v_{\sigma(m)}$ has outdegree 0, there is no edge of the form $(v_{\sigma(m)},v_{\sigma(j)})$ for $j \leq m.$  This completes the proof.
\end{proof}

\bc       \label{CSS:corollary:acyclic}
Let $G=(V,E)$ be a finite directed graph with vertex set $\{v_0,v_1,\ldots, v_{m-1} \},$ and let $\sigma$ be a permutation of  $\{0,1,\ldots,m-1 \}.$    Let $B_{\sigma}$ be the set of edges $(v_{\sigma(r)},v_{\sigma(s)}) \in E$ with $r \geq s.$  Then the subgraph $G' = (V,E\setminus B_{\sigma})$ is acyclic.
\ec

\begin{proof}
This follows immediately from Theorem~\ref{CSS:theorem:permutation}.
\end{proof}

Let $\beta(G)$ denote the minimum size of a set $X$ of edges such that the graph $G' = (V,E\setminus X)$ is directed acyclic.

\bc      \label{CSS:corollary:minacyclic}
Let $G=(V,E)$ be a finite directed graph with vertex set $\{v_0,v_1,\ldots, v_{m-1}\},$ and let $\Sigma_m$ be a set of permutations of  $\{0,1,\ldots,m-1\}.$    For $\sigma \in \Sigma_m$, let $B_{\sigma}$ be the set of edges $(v_{\sigma(r)},v_{\sigma(s)}) \in E$ with $\sigma(r)\geq \sigma(s).$  Then
\[
\beta(G) \leq \min\left\{ \card(B_{\sigma}) : \sigma \in \Sigma_m \right\}.
\]
\ec

\begin{proof}
This follows immediately from Corollary~\ref{CSS:corollary:acyclic}.
\end{proof}

Let $\gamma(G)$ denote the number of pairs of nonadjacent vertices in the undirected graph obtained from $G$ by replacing each directed edge with an undirected edge.  A \emph{tournament} is a directed graph exactly one edge between every two vertices.  If $G$ is a tournament, then $\gamma(G) = 0.$  Let $G$ be a finite, triangle-free tournament.  If $G$ contains directed cycles, then the minimum length $n$ of a directed cycle in $G$ is 4.  Let $v_{i_0}, v_{i_1},v_{i_2},\ldots, v_{i_n}$ be a cycle in $G$ of minimum length $n$.  Since $\gamma(G)=0,$ it follows that either $(v_{i_0},v_{i_2})$ or $(v_{i_2},v_{i_0})$ is an edge.  If $(v_{i_0},v_{i_2})$ is an edge, then 
$v_{i_0},v_{i_2},\ldots, v_{i_n}$ is a cycle in $G$ of length $n-1,$
which contradicts the minimality of $n$.  If $(v_{i_2},v_{i_0})$ is an edge, then $v_{i_0},v_{i_1}, v_{i_2}$ is a triangle in $G,$ which is impossible.  It follows that every tournament is directed acyclic.  Equivalently, if $G$ is triangle-free and $\gamma(G)=0,$ then $\beta(G)=0.$ 

This is a special case of a theorem of Chudnovsky, Seymour, and Sullivan\cite{chud-seym-sull07}, who proved that if $G$ is a triangle-free digraph, then $\beta(G) \leq \gamma(G).$  They conjectured that if $G$ is a 3-free digraph, then $\beta(G) \leq\gamma(G)/2.$

We shall consider the special case of the CSS conjecture in which the triangle-free graph is a Cayley graph $G=(\F_p,E_A)$ whose vertex set is the additive group of the finite field $\mathbf{F}_p$ and whose edge set $E_A$ is determined by a nonempty subset $A$ of $\mathbf{F}_p^{\ast}$ by the following rule: 
\[
E_A =  \{ ( x,x+a): x \in \mathbf{F}_p \text{ and } a \in A\}.
\]
Let $d = \card(A).$  If the Cayley graph has neither loops nor digons, then the number of pairs of adjacent vertices is the same as the number of directed edges, which is $dp$, and so the number of pairs of nonadjacent vertices is 
\[
\gamma(G) = {p \choose 2} - dp =  \frac{p(p-1-2d)}{2}.
\]
In this case the CSS conjecture asserts that 
\[
\beta(G) \leq \frac{p(p-1-2d)}{4}.
\]

\bl  \label{CSS:lemma:sigma}
Let $p$ be a prime number and $A = \{a_1,a_2,\ldots,a_d\} \subseteq \mathbf{F}_p^{\ast}$.  Let $G = (\F_p,E_A)$ be the Cayley graph constructed from $A$.  Let $\Sigma_p$ be a set of permutations of $\{0,1,2,\ldots, p-1\}$.   For $i \in  \{0,1,\ldots, p-1\}$  
 and $j \in \{1,\ldots,d\}$, define $t_{i,j} \in \{0,1,\ldots,p-1\}$ by 
\[
(\sigma(i) + p\Z )+a_j = \sigma(t_{i,j})+p\Z.
\]
Then 
\[
E_A=\{ (\sigma(i)+p\Z,\sigma(t_{i,j}) +p\Z) : i=0,\ldots,p-1 \text{ and } j=1,\ldots,d\}.
\]
Let 
\[
B_{\sigma} = \{  (\sigma(i+p\Z),\sigma(t_{i,j}+p\Z)) : t_{i,j} < i \}.
\]
The graph $G' = (\F_p, E_A \setminus B_{\sigma})$ is directed acyclic for every permutation $\sigma\in \Sigma_p$, and  
\[
\beta(G) \leq \min\{ \card(B_{\sigma}) : \sigma \in \Sigma_p \}.
\]
\el

\begin{proof}
This follows immediately from Corollary~\ref{CSS:corollary:minacyclic}.
\end{proof}

\bt  \label{CSS:theorem:beta}
Let $p$ be prime and $A = \{a_1,a_2,\ldots,a_d\} \subseteq \mathbf{F}_p^{\ast}$.  Let $G = (\F_p,E_A)$ be the Cayley graph constructed from $A$.  
Then
\[
\beta(G) \leq h_p(\langle a_1,a_2,\ldots,a_d \rangle) \leq \frac{dp}{2}.
\]
\et

\begin{proof}
Let $\Sigma_p = \{ \sigma_k\}_{k=1}^{p-1}$ be the set of permutations  of $\{0,1,2,\ldots, p-1\}$ defined by 
\[
\sigma_k(i) \equiv ki \pmod{p} \qquad\text{ for $i=0,1,\ldots,p-1.$ }
\] 
Fix $k \in \{1,2,\ldots, p-1\}$.  For $i \in  \{0,1,\ldots, p-1\}$  and $j \in \{1,\ldots,d\}$, define $t_{i,j} \in \{0,1,\ldots,p-1\} \setminus \{i\}$ by 
\[
(\sigma_k(i) + p\Z )+a_j = \sigma_k(t_{i,j})+p\Z.
\]
Let $u_k$ denote the least nonnegative integer such that $ku_k \equiv 1\pmod{p}$.  Then $\{u_1,u_2,\ldots, u_{p-1}\} = \{ 1, 2, \ldots, p-1\}.$   
Defining
\[
r_j = u_k a_j \mod p,
\]
we have $r_j \in \{1,2,\ldots, p-1\}$ and
\[
a_j = kr_j+p\Z.
\]
Then
\begin{align*}
\sigma_k(t_{i,j}) +p\Z & = (\sigma_k(i)+p\Z) + a_j \\
&  =  (ki+p\Z) + (kr_j  +p\Z) \\
& = k(i+r_j)  +p\Z \\
& = \sigma_k(i+r_j)  +p\Z
\end{align*}
and so
\[
t_{i,j} \equiv i+r_j \pmod{p}.
\]
If $i+r_j \leq p-1,$ then $t_{i,j} = i+r_j > i.$ 
If $i+r_j \geq p,$ then $t_{i,j} = i+r_j  - p <  i.$ 
It follows that $t_{i,j} < i$ if and only if $i+r_j \geq p$, that is, $p- r_j \leq i \leq p-1$
and so 
\[
\card(B_{\sigma_k}) = \sum_{j=1}^{d} r_j 
= \sum_{j=1}^{d} (u_k a_j\mod p ).
\]
By Corollary~\ref{CSS:corollary:minacyclic}, 
\begin{align*}
\beta(G) & \leq \min\{ \card(B_{\sigma_k}) : k=1,\ldots,p-1\}  \\
& = \min\left\{ \sum_{j=1}^{d} (u_k a_j\mod p )   : k=1,\ldots,p-1\right\} \\
& =  \min\left\{ \sum_{j=1}^{d} (k a_j\mod p )   : k=1,\ldots,p-1\right\} \\
& = h_p(\langle a_1,\ldots, a_d \rangle).
\end{align*}
The upper bound for the height comes from Lemma~\ref{CSS:lemma:HA}.  
\end{proof}

We  return to the CSS conjecture.  Since $dp/2 \leq p(p-1-2d)/4$
if and only if $d \leq (p-1)/4$, it follows that,  for a fixed prime $p$, we only need to consider sets $A$ of cardinality $d > p/4.$
In the other direction, Hamidoune~\cite{hami81a,nath07i} proved the Caccetta-Haggkvist conjecture for Cayley graphs:  If $A \subseteq \mathbf{F}_p^{\ast}$ and $d = |A| \geq p/r$, then the Cayley graph $(\mathbf{F}_p, E_A)$ contains a cycle of length no greater than $r$.  In particular, if the graph has no directed loops, digons, or triangles, then $d < p/3.$  Therefore, to prove the CSS conjecture for the group $\F_p$, it suffices to consider only sets $A$ of size $d$, where $p/4 < d < p/3.$

\bt
Let $p$ be a prime number, $p \geq 7,$ and let $A = \{a_1,a_2\} \subseteq \mathbf{F}_p^{\ast}$ with $a_1 \neq a_2.$   Let $G = (\F_p,E_A)$ be the Cayley graph constructed from $A$.  If $G$ is a triangle-free digraph, then 
\[
\beta(G) \leq \frac{p-1}{2} \leq \frac{\gamma(G)}{2}.
\]
\et

\begin{proof}
Since $\langle a_1,a_2 \rangle = \langle 1,a \rangle$ in  
$\PP^1(\F_p)$ with $a =a_1^{-1}a_2 \neq 1,$ and since $
\beta(G) \leq h_p(\langle a_1,a_2 \rangle) = h_p( \langle 1,a \rangle),$ it suffices to consider the case $A = \{ 1,a\}.$  
The Cayley graph $G$ is triangle-free if and only if none of the equations
\begin{align*}
x& = 0 \\
x+y & = 0\\
x+y+z& = 0
\end{align*}
has a solution with $x,y,z\in \{1,a\}.$   
The first equation implies that $a \neq 0,$ the second that $a \neq p-1+p\Z$, and that third that $2a+1 \neq 0$ and $a+2 \neq 0$, or, equivalently, that $a \neq (p-1)/2+p\Z$ or $p-2.$  It follows from Theorem~\ref{CSS:theorem:(p-1)/2} that
\[
\beta(G) \leq h_p(\langle 1,a\rangle) \leq \frac{p-1}{2} \leq \frac{p(p-5)}{4} = \frac{\gamma(G)}{2}
\]
if $p \geq 7.$  This completes the proof.
\end{proof}

\def\cprime{$'$} \def\cprime{$'$} \def\cprime{$'$}
\providecommand{\bysame}{\leavevmode\hbox to3em{\hrulefill}\thinspace}
\providecommand{\MR}{\relax\ifhmode\unskip\space\fi MR }
\providecommand{\MRhref}[2]{%
  \href{http://www.ams.org/mathscinet-getitem?mr=#1}{#2}
}
\providecommand{\href}[2]{#2}

\end{document}